\documentclass[fleqn]{mat01}
\usepackage{times,mathtimy,amssymb,latexsym}

\def\natural{\mathbf{N}}
\def\real{\mathbf{R}}
\def\complex{\mathbf{C}}

\DeclareMathOperator*{\indlim}{ind\,lim}
\DeclareMathOperator{\Spec}{Spec\,}
\DeclareMathOperator{\Dom}{Dom\,}
\DeclareMathOperator{\RE}{Re\,}

 \def\pnorm#1#2{\left|\,#2\,\right|_{#1}}
 
 \def\pNorm#1#2{\left\|\,#2\,\right\|_{#1}}
 \def\pNORM#1#2{\left|\!\left|\!\left|\,#2\,\right|\!\right|\!\right|_{#1}}

 \def\pNORMM#1#2{|\!|\!|\,#2\,|\!|\!|_{#1}}

\def\indlim{\mathop{\rm ind\,lim}}
\def\projlim{\mathop{\rm proj\,lim}}

\def\bilin#1#2{\left\langle#1,\,#2\right\rangle}
\def\Bilin#1#2{\left\langle\!\left\langle#1,\,#2
               \right\rangle\!\right\rangle}

\def\Bilinn#1#2{\langle\!\langle#1,\,#2
               \rangle\!\rangle}
\def\Biliin#1#2{\Big\langle\!\Big\langle#1,\,#2
               \Big\rangle\!\Big\rangle}

\begin{document}

\setcounter{page}{489} \firstpage{489}

\font\xxxx=msam10 at 10pt
\def\blacksquare{\mbox{\xxxx{\char'245\,}}}

\renewcommand\theequation{\thesection\arabic{equation}}

\newtheorem{theore}{Theorem}
\renewcommand\thetheore{\arabic{section}.\arabic{theore}}
\newtheorem{definit}[theore]{\rm DEFINITION}
\newtheorem{theor}[theore]{\bf Theorem}
\newtheorem{propo}[theore]{\rm PROPOSITION}
\newtheorem{lem}[theore]{\it Lemma}
\newtheorem{rem}[theore]{\it Remark}

\newtheorem{coro}[theore]{\rm COROLLARY}
\newtheorem{probl}[theore]{\it Problem}
\newtheorem{exampl}[theore]{\it Example}
\newtheorem{pot}[theore]{\it Proof of Theorem}

\def\nota{\trivlist \item[\hskip \labelsep{\it Notations.}]}

\title{Stochastic integral representations of quantum
martingales on multiple Fock space}

\markboth{Un Cig Ji}{Stochastic integral representations of
quantum martingales}

\author{UN CIG JI}

\address{Department of Mathematics, Research Institute of
Mathematical Finance, Chungbuk National University,
Cheongju~361-763,~Korea\\
\noindent E-mail: uncigji@cbucc.chungbuk.ac.kr\\[1.2pc]
\noindent {\it Dedicated to Professor Kalyan B Sinha on the
occasion of his 60th \vspace{-1pc}birthday}}

\volume{116}

\mon{November}

\parts{4}

\pubyear{2006}

\Date{}

\begin{abstract}
In this paper a quantum stochastic integral representation theorem
is obtained for unbounded regular martingales with respect to
multidimensional quantum noise. This simultaneously extends
results of Parthasarathy and Sinha to unbounded martingales and
those of the author to multidimensions.
\end{abstract}

\keyword{Fock space; quantum stochastic process; quantum
stochastic integral; quantum martingale.}

\maketitle

\section{Introduction}

The stochastic integral representations of quantum martingales
have been studied by many authors (see
\cite{BWS82,HL85,HLP86,Ji-Sinha05,Lindsay86,Meyer93,Meyer94,PS86,PS88},
etc). In \cite{PS86}, Parthasarathy and Sinha established a
stochastic integral representation of a regular bounded quantum
martingale on Fock space with respect to the basic martingales,
namely the annihilation, creation and conservation processes.
A~new proof of the Parthasarathy and Sinha representation theorem
has been discussed in \cite{Meyer94} with the special form of the
coefficient of the conservation process. In \cite{Ji03}, by using
the framework of Gaussian (white noise) analysis (see
\cite{Hida75,Obata94}), the author extended the Hudson and
Parthasarathy quantum stochastic calculus and generalized the
notion of regular martingale in the context of a certain triple of
weights \cite{Be91,LinPar89} and then the integral representation
theorem for a regular (unbounded) quantum martingale was\break
proved.

In this paper, we extend the results obtained in \cite{Ji03,PS86}
for the representation of a regular martingale to the case of
multiple Fock space with an initial Hilbert space. For our
purpose, we first extend the quantum stochastic integral studied
in \cite{HP84,PS88} (see also \cite{Par92}) to our setting.

The paper is organized as follows. In \S2 we construct a rigging
of multiple Fock space and briefly recall the basic quantum
stochastic processes. In \S3 we extend the quantum stochastic
integral studied in \cite{PS88} to a wider class of adapted
quantum stochastic processes in our setting. In \S4 we prove the
main result (Theorem~\ref{thm:main theorem}) for a stochastic
integral representation of a (unbounded) regular quantum
martingale on multiple Fock\break space.

We expect that the integral representation of quantum martingales
have applications in Markovian cocycles
\cite{AF83,GosLinSinWil03,LinWills00}. Further study is now in
progress.

\begin{nota} Let $\mathfrak{X}$ and $\mathfrak{Y}$ be
locally convex spaces.

\noindent
$\mathfrak{X}\otimes\mathfrak{Y}$: the Hilbert space tensor
product when $\mathfrak{X}$ and $\mathfrak{Y}$ are Hilbert
spaces.

\noindent $L(\mathcal{D},\mathfrak{X})$: the space of all linear
operators in $\mathfrak{X}$ with domain $\mathcal{D}$.

\noindent $\mathcal{L}(\mathfrak{X},\mathfrak{Y})$: the space of
continuous linear operators from $\mathfrak{X}$ into
$\mathfrak{Y}$ equipped with the topology of bounded convergence,
see \cite{Obata94}.
\end{nota}

\section{Multiple Fock space and basic processes}

\setcounter{equation}{0}

\setcounter{theore}{0}

Let $H=L^2(\real_+,K)\cong L^2(\real_+)\otimes K$ be the Hilbert
space of $K$-valued square integrable functions on $\real_+$ and
$B$ a selfadjoint operator in $K$ with dense domain $\Dom(B)$
satisfying $\inf\Spec(B)\ge1$, where $\real_+=[0,\infty)$ and $K$
is a separable Hilbert space called the \textit{multiplicity
space}. In fact, we take $B$ of the form
$\sum_{i\ge1}\rho_i|e_i\rangle\langle e_i|$, where $\{e_i\}$ is an
orthonormal basis for $K$ and $\{\rho_i\}$ a sequence of real
numbers greater than or equal to 1.

For each $p\in\real_+$, put
\begin{equation*}
H_p=\Dom(I\otimes B^p) \subset H
\end{equation*}
and let $H_{-p}$ be the completion of $H$ with respect to the norm
$|{I\otimes B^{-p}\cdot}|_{0}$, where $\pnorm{0}{\cdot}$ is the
norm on $H$. Then we have
\begin{equation*}
H_\infty =\projlim_{p\rightarrow\infty}H_p\subset H\cong
H^*\subset H_\infty^*\cong H_{-\infty}
=\indlim_{p\rightarrow\infty}H_{-p},
\end{equation*}
where $H_\infty^*$ is the strong dual space of $H_\infty$ with
respect to $H$.

The (Boson) Fock space over $H$ is denoted by
$\mathcal{H}=\Gamma(H)$. Then by definition, $\mathcal{H}$ is the
space of sequences $\phi=(f_n)_{n=0}^\infty$, where $f_n\in
H^{\widehat\otimes n}$ ($n$-fold symmetric tensor power of the
Hilbert space $H$) such that
\begin{equation*}
\pNorm{0}{\phi}^2 =\sum_{n=0}^\infty n!|{f_n}|_{0}^2<\infty,
\end{equation*}
where $|\cdot|_0$ is the norm on $H^{\widehat\otimes n}$ for any
$n\in\natural$.

Let $\mathcal{I}$ be a separable Hilbert space which is called the
\textit{initial Hilbert space} and $A$ a selfadjoint operator in
$\mathcal{I}$ with dense domain $\Dom(A)$ satisfying
$\inf\Spec(A)\ge1$. To lighten the notation, the operator
$A\otimes\Gamma(eI\otimes B)$ in $\mathcal{I}\otimes\mathcal{H}$
is denoted by $\mathbf{A}$ and
\begin{equation*}
\mathbf{A}^\mathbf{p}=A^{p_1}\otimes\Gamma(e^{p_2}I\otimes
B^{p_3}),\quad \mathbf{p}=(p_1,p_2,p_3)\in\real^3,
\end{equation*}
where $\Gamma(C)$ is the second quantization of the operator $C$
(see \cite{Par92}). Then by standard arguments we may construct a
triplet:
\begin{equation*}
\mathcal{G}_\infty\subset\mathcal{G}\subset\mathcal{G}_{-\infty}
\end{equation*}
from $\mathcal{G}=\mathcal{I}\otimes\mathcal{H}$ and
$\mathbf{A}=A\otimes\Gamma(eI\otimes B)$. More precisely, for each
$\mathbf{p}\in\real_+^3$, put
\begin{equation*}
\mathcal{G}_\mathbf{p}=\Dom(\mathbf{A}^\mathbf{p})
\subset\mathcal{G}\equiv\mathcal{I}\otimes\mathcal{H}
\end{equation*}
and then $\mathcal{G}_{\mathbf{p}}$ becomes a Hilbert space with
norm
$\pNORM{\mathbf{p}}{\cdot}=\pNORM{0}{\mathbf{A}^{\mathbf{p}}\cdot}$,
where $\pNORM{0}{\cdot}$ is the norm on
$\mathcal{I}\otimes\mathcal{H}$. Let $\mathcal{G}_{-\mathbf{p}}$
be the completion of $\mathcal{I}\otimes\mathcal{H}$ with respect
to the norm
$\pNORM{\mathbf{-p}}{\cdot}=|\!|\!|{\mathbf{A}^{-\mathbf{p}}\cdot}|\!|\!|_{0}$,
and
\begin{equation*}
\mathcal{G}_\infty
=\projlim_{p_1,p_2,p_3\rightarrow\infty}\mathcal{G}_{(p_1,p_2,p_3)},
\quad \mathcal{G}_{-\infty}
=\indlim_{p_1,p_2,p_3\rightarrow\infty}\mathcal{G}_{(-p_1,-p_2,-p_3)}.
\end{equation*}
Note that $\mathcal{G}_{-\infty}$ is topologically isomorphic to
the strong dual space $\mathcal{G}_\infty^*$ of
$\mathcal{G}_\infty$ with respect to
$\mathcal{I}\otimes\mathcal{H}$.

For each $\mathbf{p}=(p_2,p_3)\in\real_+^2$, put
\begin{equation*}
\mathcal{H}_\mathbf{p}=\Dom(\Gamma(e^{p_2}I\otimes B^{p_3}))
\end{equation*}
and let $\mathcal{H}_{-\mathbf{p}}$ be the completion of
$\mathcal{H}$ with respect to the norm
$\pNorm{\mathbf{-p}}{\cdot}=\|\Gamma(e^{-p_2}I\otimes$\break
${B}^{-p_3})\cdot\|_{0}$, and
\begin{equation*}
\mathcal{H}_\infty
=\projlim_{p_2,p_3\rightarrow\infty}\mathcal{H}_{(p_2,p_3)}, \quad
\mathcal{H}_{-\infty}
=\indlim_{p_2,p_3\rightarrow\infty}\mathcal{H}_{(-p_2,-p_3)}.
\end{equation*}

For each interval $[a,b]\subset\real_+$, we write
$H_{[a,b]}=L^2([a,b],K)$ and then
\begin{equation*}
H=H_{s]}\oplus H_{[s,t]}\oplus H_{[t}, \quad 0<s<t<\infty
\end{equation*}
with abbreviations $H_{s]}$ and $H_{[t}$ when $[a,b]=[0,s]$ and
$[a,b]=[t,\infty]$, respectively. Therefore, we have the
identification
\begin{equation*}
\mathcal{G}
=\mathcal{G}_{s]}\otimes\mathcal{H}_{[s,t]}\otimes\mathcal{H}_{[t},
\quad \mathcal{G}_{s]}=\mathcal{I}\otimes\mathcal{H}_{s]},
\end{equation*}
where
\begin{equation*}
\mathcal{H}_{s]}={\Gamma}(H_{s]}), \quad
\mathcal{H}_{[s,t]}={\Gamma}(H_{[s,t]}),\quad
\mathcal{H}_{[t}={\Gamma}(H_{[t}).
\end{equation*}
Moreover, for any
$\mathbf{p}=(p_1,p_2,p_3)\in\real_+^3\cup\real_-^3$
($\real_-=(-\infty,0]$) and $0<s<t<\infty$, we have
\begin{equation*}
\mathcal{G}_\mathbf{p}
=\mathcal{G}_{\mathbf{p};s]}\otimes\mathcal{H}_{\mathbf{p}';[s,t]}
\otimes\mathcal{H}_{\mathbf{p}';[t},
\end{equation*}
where $\mathbf{p}'=(p_2,p_3)$ and
\begin{equation*}
\mathcal{G}_{\mathbf{p};s]}=\mathcal{G}_\mathbf{p}\cap\mathcal{G}_{s]},
\quad\mathcal{H}_{\mathbf{p}';[s,t]}=\mathcal{H}_{\mathbf{p}'}\cap\mathcal{H}_{[s,t]},
\quad\mathcal{H}_{\mathbf{p}';[t}=\mathcal{H}_{\mathbf{p}'}\cap\mathcal{H}_{[t}
\end{equation*}
(closures when $\mathbf{p}\in\real_-^3$).

For each $g,h\in H_\infty$ and
$T\in\mathcal{L}(H_\infty,H_\infty)$, the \textit{annihilation},
\textit{creation} and \textit{conservation operators} are defined
on $\mathcal{H}_\infty$ as follows:
\begin{align*}
a(g)\phi &= (ng\widehat{\otimes}^1 f_n)_{n=1}^\infty,\nonumber\\[.4pc]
a^*(h)\phi &= (S_{1+n}(h\otimes f_n))_{n=0}^\infty,\nonumber\\[.4pc]
\lambda(T)\phi &= ((n+1)S_{1+n}(T\otimes I^{\otimes n})f_{n+1}
)_{n=0}^\infty,\nonumber
\end{align*}
respectively, for any
$\phi=(f_n)_{n=0}^\infty\in\mathcal{H}_\infty$, where
$g\widehat{\otimes}^1 f_n$ is the left $1$-contraction of $g$ and
$f_n$ \cite{Obata94}, and $S_{l+m}$ stands for the symmetrizing
operator. Then we can easily show that $a(g)$, $a^*(h)$ and
$\lambda(T)$ are continuous linear operators acting on
$\mathcal{H}_\infty$. The operators $a(g)$ and $a^*(g)$ are
adjoint to each other  and $\lambda(T^*)=(\lambda(T))^*$.

The three basic (quantum stochastic) processes called
\textit{annihilation}, \textit{creation} and \textit{conservation
processes} are defined by
\begin{align*}
A_i(t)&= I\otimes a(\mathbf{1}_{[0,t]}\otimes e_i), \nonumber\\[.4pc]
A_i^*(t)&= I\otimes a^*(\mathbf{1}_{[0,t]}\otimes e_i),\nonumber\\[.4pc]
\Lambda_{ij}(t)&= I\otimes\lambda(\mathbf{1}_{[0,t]}\otimes
P_{ij}), \nonumber
\end{align*}
respectively, where $I$ is the identity operator on $\mathcal{I}$
and $\mathbf{1}_{[0,t]}$ the indicator function. In the definition
of $A_i(t)$ and $A_i^*(t)$ the indicator function
$\mathbf{1}_{[0,t]}$ is a vector in $L^2(\real_{+})$ while it is
considered as a multiplication operator in $L^2(\real_{+})$ in the
definition of $\Lambda_{ij}(t)$.

For each $f\in H$, a vector of the form:
\begin{equation*}
\phi_f=\left(1,f,\frac{f^{\otimes 2}}{2!},\dots, \frac{f^{\otimes
n}}{n!}, \dots\right)
\end{equation*}
is called an \textit{exponential vector} or a \textit{coherent
vector}. Note that $\phi_f$ belongs to $\mathcal{H}_\infty$ (resp.
$\mathcal{H}_{-\infty}$) if and only if $f$ belongs to $H_\infty$
(resp. $H_{-\infty}$). The exponential vectors $\{\phi_f\,;\,f\in
H_\infty\}$ span a dense subspace of $\mathcal{H}_\infty$, hence
of $\mathcal{H}_\mathbf{p}$ for all $\mathbf{p}\in\real_+^2$ and
of $\mathcal{H}_{-\infty}$. We denote $\mathcal{E}(D)$ the linear
subspace generated by $\{\phi_f\,;\, f\in D\}$ for $D\subset H$.
Then for any $f,g\in H_\infty$ and $t\in\real_+$ we have
\begin{align*}
\Bilinn{A_i(t)u\otimes\phi_f}{v\otimes\phi_g}
&=\bilin{u}{v}\left(\int_0^tf_i(s)\hbox{d}s\right)e^{\bilin{f}{g}},\\[.3pc]
\Bilinn{A_i^*(t)u\otimes\phi_f}{v\otimes\phi_g}
&=\bilin{u}{v}\left(\int_0^tg_i(s)\hbox{d}s\right)e^{\bilin{f}{g}},\\[.3pc]
\Bilinn{\Lambda_{ij}(t)u\otimes\phi_f}{v\otimes\phi_g}
&=\bilin{u}{v}\left(\int_0^tf_j(s)g_i(s)\hbox{d}s\right)e^{\bilin{f}{g}},
\end{align*}
where $\Bilin{\cdot}{\cdot}$ is the $\complex$-bilinear form on
$\mathcal{G}_{-\infty}\times\mathcal{G}_\infty$ and
$h_i(s)=\bilin{h(s)}{e_i}$ for $h\in H$.

The quantum Ito's formula established in \cite{HP84} is summarized
by the following table:
\renewcommand{\arraystretch}{1.3}
\begin{equation}\label{eqn:quantum Ito formula}
\begin{tabular}{|c|c|c|c|c|}\hline
 {}                &$\hbox{d}t$&$\hbox{d}A_k$&$\hbox{d}A_k^{*}$            &$\hbox{d}\Lambda_{kl}$             \\ \hline
$\hbox{d}t$               & 0  & 0    & 0                    & 0                          \\ \hline
$\hbox{d}A_{i}$           & 0  & 0    & $\delta_{ik}\hbox{d}t$      & $\delta_{ik}\hbox{d}A_{l}$        \\ \hline
$\hbox{d}A_{i}^{*}$       & 0  & 0    & 0                    & 0                          \\ \hline
$\hbox{d}\Lambda_{ij}$    & 0  & 0    & $\delta_{jk}\hbox{d}A_i^{*}$& $\delta_{jk}\hbox{d}\Lambda_{il}$ \\ \hline
\end{tabular}\ \ .
\end{equation}

\section{Quantum stochastic integral}

\setcounter{equation}{0}

\setcounter{theore}{0}

Let $\mathcal{D}_0$ and $M$ be dense linear subspaces of
$\mathcal{I}_\infty$ and $H_\infty$, respectively, such that
$\mathbf{1}_{[0,t]}f\in M$ for any $t\in\real_+$ and $f\in M$, and
let
\begin{equation*}
M_{t]}=\{\mathbf{1}_{[0,t]}f\,;\,f\in M\}, \quad
M_{[t}=\{\mathbf{1}_{[t,\infty)}f\,;\,f\in M\}.
\end{equation*}
We put $\widetilde{\mathcal{E}}=\mathcal{D}_0\otimes_{\rm
al}\mathcal{E}(M)\subset\mathcal{G}_\infty$, where $\otimes_{\rm
al}$ is the algebraic tensor product, and put
\begin{equation*}
\widetilde{\mathcal{E}}_{t]}=\mathcal{D}_0\otimes_{\rm
al}\mathcal{E}(M_{t]}),
\quad \mathcal{E}_{[t}=\mathcal{E}(M_{[t})\quad\text{and then}\quad
\widetilde{\mathcal{E}}
=\widetilde{\mathcal{E}}_{t]}\otimes_{\rm al}\mathcal{E}_{[t}.
\end{equation*}

A family of operators $\Xi=\{\Xi(t)\}_{t\ge0}\subset
L(\widetilde{\mathcal{E}},\mathcal{G}_{-\infty})$ is called a
\textit{$\mathcal{G}_\mathbf{p}$-quantum stochastic process} if
there exists $\mathbf{p}\in\real_+^3\cup\real_-^3$ (independent of
$t\ge0$) such that $\Xi(t)\in
L(\widetilde{\mathcal{E}},\mathcal{G}_\mathbf{p})$ for each
$t\ge0$ and for each $\psi\in\widetilde{\mathcal{E}}$ the map
$\real_+\ni t\mapsto\Xi(t)\psi\in\mathcal{G}_\mathbf{p}$ is
strongly measurable. We may then think of $\Xi(t)$ as a densely
defined operator on the Hilbert space $\mathcal{G}_\mathbf{p}$;
and call $\Xi$ \textit{adapted} if $\Xi(t)=\Xi(t])\otimes_{\rm
alg}I([t)$ for some $\Xi(t])\in
L(\widetilde{\mathcal{E}}_{t]},\mathcal{G}_{\mathbf{p};t]})$,
where $I([t)$ is the identity operator on
$\mathcal{G}_{\mathbf{p};[t}$.

For certain sets of $\{E^{(k)}\}_{k=1,2,3,4}$ of families of
adapted process, stochastic integrals of the type
\begin{equation*}
\int_0^t\left\{\sum_{i,j} E_{ij}^{(1)}\hbox{d}\Lambda_{ij} +\sum_{i}
E_i^{(2)}\hbox{d}A_i +\sum_{i} E_i^{(3)}\hbox{d}A_{i}^{+} +E^{(4)}\hbox{d}s\right\}
\end{equation*}
can be defined as in \cite{Ji03}. We first define the integrals
for a finite family of simple adapted processes
$\{E^{(k)}\}_{k=1,2,3,4}$ and then the definition can be extended
to a certain class of countable families $\{E^{(k)}\}_{k=1,2,3,4}$
with a norm estimate (see \eqref{eqn:norm estimate of QSI})
induced by the quantum It\^o formula. For detailed calculations, we
refer to \cite{HP84} and \cite{PS88}.

For each $\mathbf{p}=(p_1,p_2,p_3)\in\real_+^3\cup\real_-^3$ we denote
$\mathcal{A}_2(\widetilde{\mathcal{E}},\mathcal{G}_\mathbf{p})$
the class of all (ordered) quadruples of families of adapted
processes
\begin{equation*}
\mathbf{E}\equiv\{E_{ij}^{(1)}(t),\,\,E_i^{(2)}(t),\,\,
E_i^{(3)}(t),\,\,E^{(4)}(t);\,\, 1 \le i,j<\infty,\,\,t\ge0\}
\end{equation*}
satisfying
\begin{align}\label{eqn:condition for QSI}
&\int_0^t \left\{\sum_{i}\rho_i^{2p_3}
\pNORM{\mathbf{p}}{\sum_{j}f_j(s)E_{ij}^{(1)}(s)u\otimes\phi_f}^2\right.\nonumber\\[.3pc]
&\quad\,\left. +\sum_{k=2}^3\sum_{i}\rho_i^{2p_3}
\pNORMM{\mathbf{p}}{E_i^{(k)}(s)u\otimes\phi_f}^2
+\pNORMM{\mathbf{p}}{E^{(4)}(s)u\otimes\phi_f}^2\right\}\hbox{d}s<\infty
\end{align}
for all $t>0$, $u\in\mathcal{D}_0$ and $f\in M$.

\begin{theor}[\!]\label{thm:QSI}
Let $\mathbf{p}=(p_1,p_2,p_3)\in\real_+^3\cup\real_-^3$ and
$\mathbf{E}\in\mathcal{A}_2(\widetilde{\mathcal{E}},\mathcal{G}_\mathbf{p})$.
Then the stochastic integral
\begin{align*}
\Xi(t) &= \int_0^t\sum_{i,j} E_{ij}^{(1)} (s){\rm d}\Lambda_{ij}(s)
+\int_0^t\sum_{i} E_i^{(2)} (s){\rm d}A_i(s)\nonumber\\[.3pc]
&\quad\, +\int_0^t\sum_{i}
E_i^{(3)}(s){\rm d}A_{i}^{+}(s)+\int_0^tE^{(4)}(s){\rm d}s\nonumber
\end{align*}
is well-defined as an adapted process in
$L(\widetilde{\mathcal{E}},\mathcal{G}_\mathbf{p})$. Moreover{\rm
,} for any $u,v\in\mathcal{D}_0$ and $f,g\in M$ we have
\begin{align}
&\Bilinn{\Xi(t)u\otimes\phi_f}{v\otimes\phi_g}\nonumber\\[.3pc]
&\ \ =\int_0^t \Bilin{\left\{\sum_{i,j}
g_i(s)f_j(s)E_{ij}^{(1)}(s)
+\!\sum_{i} f_i(s)E_i^{(2)}(s)\right\}u\!\otimes\phi_f}{v\otimes\phi_g}{\rm d}s\nonumber\\[.3pc]
&\quad\, +\int_0^t \Bilin{\left\{\sum_{i} g_i(s)E_i^{(3)}(s)
+E^{(4)}(s)\right\}u\otimes\phi_f}{v\otimes\phi_g}{\rm d}s
\label{eqn:Bilinear form of QSI}
\end{align}
and
\begin{equation}\label{eqn:norm estimate of QSI}
\pNORMM{\mathbf{p}}{\Xi(t)u\otimes\phi_f}^2
\le\exp\left\{t+3e^{2p_2}\int_0^t\pnorm{p_3}{f(u)}^2{\rm
d}u\right\} \left(\int_0^tG(s){\rm d}s\right) <\infty,
\end{equation}
where{\rm ,} for each $t\in\real_+${\rm ,}
\begin{align}
G(t) &= 3e^{2p_2}\sum_i\rho_i^{2p_3}\pNORM{\mathbf{p}}{\sum_j
f_j(t)E_{ij}^{(1)}(t)u\otimes\phi_f}^2\nonumber\\[.3pc]
&\quad\, +
e^{2p_2}\sum_i\rho_i^{2p_3}\pNORMM{\mathbf{p}}{E_i^{(2)}(t)u\otimes\phi_f}^2
\nonumber\\[.3pc]
&\quad\, +
3e^{2p_2}\sum_i\rho_i^{2p_3}\pNORMM{\mathbf{p}}{E_i^{(3)}(t)u\otimes\phi_f}^2
    +\pNORMM{\mathbf{p}}{E^{(4)}(t)u \otimes\phi_f}^2.
    \label{def:G(t)}
\end{align}
\end{theor}

\begin{proof}
By similar arguments of those used in \cite{PS88} and \cite{Ji03}
using the quantum It\^o formula \eqref{eqn:quantum Ito formula},
for simple quadruple $\mathbf{E}$ with finite number of non-zero
components we can compute that
\begin{align}\label{eqn:eqn1 for norm estimate of QSI}
\frac{\rm d}{{\rm d}t}\pNORMM{\mathbf{p}}{\Xi(t)u\otimes\phi_f}^2
&=2\RE\left\{\sum_{k=1}^5 S_k\right\}\nonumber\\[.3pc]
&\quad\, + \sum_ie^{2p_2}\rho_i^{2p_3}
   \pNORM{\mathbf{p}}{\sum_jf_j(t) E_{ij}^{(1)}(t)u \otimes\phi_f}^2
   \nonumber\\[.3pc]
&\quad\, +\sum_ie^{2p_2}\rho_i^{2p_3}
   \pNORMM{\mathbf{p}}{E_i^{(3)}(t)u \otimes\phi_f}^2,
\end{align}
where
\begin{align*}
\hskip -4pc S_1&=e^{2p_2}\sum_{i}\rho_i^{2p_3}
  \Bilin{\mathbf{A}^\mathbf{p}(f_i(t)\Xi(t)u\otimes\phi_f)}
        {\overline{\mathbf{A}^\mathbf{p}\left(\sum_{j}f_j(t)E_{ij}^{(1)}(t)u\otimes\phi_f\right)}},\\[.4pc]
\hskip -4pc S_2&=\sum_{i}e^{2p_2}\rho_i^{2p_3}
  \Biliin{\mathbf{A}^\mathbf{p}(\bar{f}_i(t)\Xi(t)u\otimes\phi_f)}
        {\overline{\mathbf{A}^\mathbf{p}(E_i^{(2)}(t)u\otimes\phi_f)}},\\[.3pc]
\hskip -4pc S_3&=\sum_{i}e^{2p_2}\rho_i^{2p_3}
  \Biliin{\mathbf{A}^\mathbf{p}(f_i(t)\Xi(t)u\otimes\phi_f)}
        {\overline{\mathbf{A}^\mathbf{p}(E_i^{(3)}(t)u\otimes\phi_f)}},\\[.3pc]
\hskip -4pc
S_4&=\Bilin{\mathbf{A}^\mathbf{p}(\Xi(t)u\otimes\phi_f)}
         {\overline{\mathbf{A}^\mathbf{p}(E^{(4)}(t)u\otimes\phi_f)}},\\[.4pc]
\hskip -4pc S_5&=\sum_{i}e^{2p_2}\rho_i^{2p_3}
  \Bilin{\mathbf{A}^\mathbf{p}(E_i^{(3)}(t)u\otimes\phi_f)}
        {\overline{\mathbf{A}^\mathbf{p}\left(\sum_jf_j(t)E_{ij}^{(1)}u\otimes\phi_f\right)}}.
\end{align*}
By using the Cauchy--Schwarz inequality and the fact
$2\RE\bar{a}b\le|a|^2+|b|^2$, we obtain from \eqref{eqn:eqn1 for
norm estimate of QSI} that
\begin{equation*}
\frac{\rm d}{{\rm d}t}\pNORMM{\mathbf{p}}{\Xi(t)u\otimes\phi_f}^2
\le
(1+3e^{2p_2}|{f(t)}|_{{p_3}}^2)\pNORMM{\mathbf{p}}{\Xi(t)u\otimes\phi_f}^2
+G(t),
\end{equation*}
where $G(t)$ is given as in \eqref{def:G(t)}. The inequality
\eqref{eqn:norm estimate of QSI} can be obtained by applying
Gronwall's lemma with the above inequality, as in \cite{HP84} or
\cite{Par92}. Then the inequality \eqref{eqn:norm estimate of QSI}
allows the extension of the integral to
$\mathcal{A}_2(\widetilde{\mathcal{E}},\mathcal{G}_\mathbf{p})$
satisfying the inequality \eqref{eqn:norm estimate of QSI}.\hfill
$\blacksquare$
\end{proof}

\section{Regular quantum martingales}

\setcounter{equation}{0} \setcounter{theore}{0}

An adapted processes $\{\Xi(t)\}_{t\ge0}\subset
L(\widetilde{\mathcal{E}},\mathcal{G}_\mathbf{p})$ is called a
\textit{quantum martingale} if for any $0\le s\le t$,
\begin{align*}
\Bilinn{\Xi(t)(u\otimes\phi_{\mathbf{1}_{[0,s]}f})}{v\otimes\phi_{\mathbf{1}_{[0,s]}g}}
=\Bilinn{\Xi(s)(u\otimes\phi_{\mathbf{1}_{[0,s]}f})}{v\otimes\phi_{\mathbf{1}_{[0,s]}g}}
\end{align*}
for any $u,v\in\mathcal{D}_0$ and $f,g\in M$. For each $1\le
i,j<\infty$, the annihilation process $\{A_i(t)\}_{t\ge0}$,
creation process $\{A_i^*(t)\}_{t\ge0}$ and conservation process
$\{\Lambda_{ij}(t)\}_{t\ge0}$ are quantum martingales which are
called the basic martingales in quantum stochastic calculus.

In the following, for
$\mathbf{p},\mathbf{q}\in\real_+^3\cup\real_-^3$ with
$\mathbf{p}-\mathbf{q}\in\real_+^3$ we consider quantum
martingales $\Xi$ in
$\mathcal{L}(\mathcal{G}_\mathbf{p},\mathcal{G}_\mathbf{q})$.
Thus, for any $0\le s\le t$ and $\phi_{s]}\in
\mathcal{G}_{\mathbf{p};s]}$, $\psi_{s]}\in
\mathcal{G}_{-\mathbf{q};s]}$,
\begin{equation*}
\Bilinn{\Xi_t\phi_{s]}}{\psi_{s]}}=\Bilinn{\Xi_s\phi_{s]}}{\psi_{s]}}.
\end{equation*}
The following definition of regular martingale is a simple
modification of the definition of bounded regular martingale in
\cite{PS86} and \cite{Ji03}.

\begin{definit}$\left.\right.$\vspace{.5pc}

\noindent {\rm A quantum martingale $\Xi$ in
$\mathcal{L}(\mathcal{G}_\mathbf{p},\mathcal{G}_\mathbf{q})$ is
said to be \textit{regular with respect to a Radon measure
$\mathfrak{m}$} on $\real_+$, or simply \textit{regular} if for
any $0\le v<u$ and $\phi\in \mathcal{G}_{\mathbf{p};v]}$,
$\psi\in\mathcal{G}_{-\mathbf{q};v]}$,
\begin{align}\label{eqn:regularity condition}
&\pNORM{\mathbf{q}}{(\Xi_{u}-\Xi_{v})\phi}^2
 \le\pNORM{\mathbf{p}}{\phi}^2\,\mathfrak{m}([v,u]),\nonumber\\[.3pc]
&\pNORMM{-\mathbf{p}}{(\Xi_{u}^*-\Xi_{v}^*)\psi}^2
 \le\pNORM{-\mathbf{q}}{\psi}^2\,\mathfrak{m}([v,u]).
\end{align}}
\end{definit}

\begin{propo}\label{pro:regularity of
martingale}$\left.\right.$\vspace{.5pc}

\noindent Let $\Xi$ be a quantum martingale in
$\mathcal{L}(\mathcal{G}_\mathbf{p},\mathcal{G}_\mathbf{q})$. If
$\Xi$ has the integral representation$:$
\begin{equation*}
\hbox{\rm d}\,\Xi =\sum_{i,j} E_{ij}\hbox{\rm
d}\Lambda_{ij}+\sum_{i}F_i^*\hbox{\rm d}A_i+\sum_{i}G_i \hbox{\rm
d}A_i^*,
\end{equation*}
where the quadruples $(E_{ij},F_i^*,G_i,0)$ and
$(E_{ij}^*,F_i,G_i^*,0)$ belong to
$\mathcal{A}_2(\widetilde{\mathcal{E}},\mathcal{G}_\mathbf{q})$
and
$\mathcal{A}_2(\widetilde{\mathcal{E}},\mathcal{G}_\mathbf{-p}),$
respectively{\rm ,} and $E_{ij},F_i^*,G_i$ are adapted processes
in $\mathcal{L}(\mathcal{G}_\mathbf{p},\mathcal{G}_\mathbf{q})$
such that
\begin{equation*}
\sum_{i}\rho_i^{2q_3}G_i^\dagger(s)\mathbf{A}^{2\mathbf{q}}G_i(s)
\quad\text{and}\quad
\sum_{i}\rho_i^{-2p_3}F_i^\dagger(s)\mathbf{A}^{-2\mathbf{p}}F_i(s)
\end{equation*}
converge weakly to self-adjoint operators
$G(s)\in\mathcal{L}(\mathcal{G}_\mathbf{p},\mathcal{G}_{-\mathbf{p}})$
and
$F(s)\in\mathcal{L}(\mathcal{G}_{-\mathbf{q}},\mathcal{G}_{\mathbf{q}}),$
respectively{\rm ,} with the property that
$\|G(s)\|_{\mathbf{p};-\mathbf{p}}$ and
$\|F(s)\|_{-\mathbf{q};\mathbf{q}}$ are locally integrable{\rm ,}
where $K^\dagger$ denotes the adjoint of the operator $K$ with
respect to $\Bilin{\cdot}{\overline{\,\cdot\,}}$ and
$\|\Xi\|_{\mathbf{r};\mathbf{s}}$ is the operator norm of
$\Xi\in\mathcal{L}(\mathcal{G}_\mathbf{r},\mathcal{G}_\mathbf{s})$.
Then $\Xi$ is regular.
\end{propo}

\begin{proof}
Let $\mathbf{p}=(p_1,p_2,p_3)\in\real_+^3\cup\real_-^3$ and
$\mathbf{q}=(q_1,q_2,q_3)\in\real_+^3\cup\real_-^3$. Note that for
any $0\le a<t$ and $\phi_{a]}\in \mathcal{G}_{\mathbf{p};a]}$,
$\psi_{a]}\in\mathcal{G}_{\mathbf{q};a]}$,
\begin{equation*}
\Bilinn{\Xi(t)\phi_{a]}}{\psi_{a]}}_\mathbf{q}
=\Bilinn{\Xi(s)\phi_{a]}}{\psi_{a]}}_\mathbf{q}.
\end{equation*}
It follows that
\begin{equation}\label{eqn:basic equality for regularity-0}
\pNORMM{\mathbf{q}}{(\Xi(t)-\Xi(a))\phi_{a]}}^2
=\pNORMM{\mathbf{q}}{\Xi(t)\phi_{a]}}^2-\pNORMM{\mathbf{q}}{\Xi(a)\phi_{a]}}^2.
\end{equation}
Therefore, by applying \eqref{eqn:eqn1 for norm estimate of QSI},
we obtain that for any $0\le a<t$ and
$\phi_{a]}\in\widetilde{\mathcal{E}}_{a]}$,
\begin{align}
\pNORMM{\mathbf{q}}{(\Xi(t)-\Xi(a))\phi_{a]}}^2
&=e^{2q_2}\int_a^t\sum_{i}\rho_i^{2q_3}
   \pNORMM{\mathbf{q}}{G_i(s)\phi_{a]}}^2{\rm d}s\nonumber\\[.3pc]
&\le e^{2q_2}\pNORMM{\mathbf{p}}{\phi_{a]}}^2
   \int_a^t \|G(s)\|_{\mathbf{p};-\mathbf{p}}{\rm d}s.
   \label{eqn:basic equality for regularity}
\end{align}
Similarly, for any $0\le a<t$ and
$\psi_{a]}\in\widetilde{\mathcal{E}}_{a]}$, we have
\begin{equation}\label{eqn:basic equality for regularity1}
\pNORMM{-\mathbf{p}}{(\Xi(t)^*-\Xi(a)^*)\psi_{a]}}^2
  \le e^{-2p_2}\pNORMM{-\mathbf{q}}{\psi_{a]}}^2
  \int_a^t \|F(s)\|_{-\mathbf{q};\mathbf{q}}^2 {\rm d}s.
\end{equation}
Now, we define a Radon measure $\mathfrak{m}$ on $\real_+$ by
\begin{align*}
\mathfrak{m}([a,b]) &= \int_a^b
(e^{2q_2}\|G(s)\|_{\mathbf{p};-\mathbf{p}}^2
+e^{-2p_2}\|F(s)\|_{-\mathbf{q};\mathbf{q}}^2){\rm d}s\\[.3pc]
&\quad\, {\rm for\,\, all\,\,} 0\le a\le b<\infty.
\end{align*}
Therefore, by (\ref{eqn:basic equality for regularity}),
 (\ref{eqn:basic equality for regularity1}) and
 the density of $\widetilde{\mathcal{E}}_{a]}$ in $\mathcal{G}_{\mathbf{p};a]}$
 and $\mathcal{G}_{-\mathbf{q};a]}$, we see that $\Xi$ is regular
 with respect to the absolutely continuous Radon measure
 $\mathfrak{m}$.\hfill $\blacksquare$
\end{proof}

\begin{rem}\label{remark:remark1 for IRT}
{\rm Let $\Xi$ be a martingale in
$\mathcal{L}(\mathcal{G}_\mathbf{p},\mathcal{G}_\mathbf{q})$ which
is regular with respect to the Radon measure $\mathfrak{m}$. Then
for any $t>a$,
\begin{align*}
\|\Xi(t)\|_{\mathbf{p};\mathbf{q}}
&\ge\sup_{\pNORM{\mathbf{p}}{\phi_{a]}}=1}
        \pNORMM{\mathbf{q}}{\Xi(t)\phi_{a]}}\\[.3pc]
&\ge \sup_{\pNORM{\mathbf{p}}{\phi_{a]}}=1}
        \pNORMM{\mathbf{q}}{\Xi(a)\phi_{a]}}
\ge\|\Xi(a)\|_{\mathbf{p};\mathbf{q}},
\end{align*}
where we used (\ref{eqn:basic equality for regularity-0}) for the
second inequality. Therefore,
$\|\Xi(\cdot)\|_{\mathbf{p};\mathbf{q}}$ is non-decreasing. }
\end{rem}

Let $P$ denote the probability measure of an independent
identically distributed sequence $\{B_1,B_2,\dots\}$ of standard
Brownian motions. Then the Hilbert space $L_2(P)$ is identified
with $\Gamma(L_2(\real_+,\real)\otimes\ell_2)$ by the following
correspondence:
\begin{equation*}
\phi_f\,\,\longleftrightarrow\,\, \exp\sum_{i}\left(\int_0^\infty
f_i\,{\rm d}B_i -\frac{1}{2}\int_0^\infty f_i^2\,{\rm d}t\right),
\end{equation*}
where ${f}=(f_1,f_2,\dots)\in\bigoplus_{i=1}^\infty
L_2(\real_+,\real)\cong L_2(\real_+,\real)\otimes\ell_2$. Put
\begin{align*}
M_0 &= \{f=(f_1,\ldots,f_i,\ldots)\in H_\infty\,;\\[.3pc]
&\quad f_i=0 \ \text{ for all but a finite number of }
i'\text{s}\}.
\end{align*}
Then $\mathcal{E}_0=\mathcal{E}(M_0)$ and
$\widetilde{\mathcal{E}}_0=\mathcal{I}_\infty\otimes_{\rm
al}\mathcal{E}(M_0)$ are total in $H$ and $\mathcal{G}$,
respectively, where $\mathcal{I}_\infty$ is the Fr\'echet space
constructed by the standard manner with $\mathcal{I}$ and the
positive operator $A$, and then we have
\begin{equation*}
\phi_{f\mathbf{1}_{[0,t]}}-1
=\sum_{i}\int_0^tf_i(s)\phi_{f\mathbf{1}_{[0,s]}}\,{\rm d}B_i(s)
\end{equation*}
for $\phi_f\in\mathcal{E}_0$. In general, we have the following
proposition which is an extension of the classical martingale
representation theorem of Kunita--Watanabe for $L^2$-martingales
adapted to one Brownian motion to an $\mathcal{I}$-valued
$L^2$-martingale adapted to a countable family of independent
Brownian motion.

\begin{propo}{\rm \cite{PS88}}\label{prop:extension of K-W
Theorem}$\left.\right.$\vspace{.5pc}

\noindent Let $\{X(t)\}_{t\ge0}$ be an $\mathcal{I}$-valued square
integrable martingale adapted to $\{B_i\}$ which is an independent
identically distributed sequence of standard Brownian motions.
Then
\begin{equation*}
X(t)=X(0)+\sum_{i}\int_0^t\xi_i\,{\rm d}B_i,
\end{equation*}
where $\{\xi_i\}_{i\ge1}$ is a sequence of adapted processes
satisfying
\begin{equation*}
\int_0^t\sum_{i}\mathbf{E}[\|\xi_i(s)\|_\mathcal{I}^2]{\rm d}s
<\infty,\quad t\in\real_+.
\end{equation*}
\end{propo}

Our aim is to prove the converse of
Proposition~\ref{pro:regularity of martingale} generalizing the
main result in \cite{Ji03} and \cite{PS88}. For the proof we use
similar arguments to those used in \cite{PS88}.

\begin{theor}[\!]\label{thm:main theorem}
Let $\mathbf{p},\mathbf{q}\in\real_+^3\cup\real_-^3$ with
$\mathbf{p}-\mathbf{q}\in\real_+^3$. Let $\Xi$ be a martingale in
$\mathcal{L}(\mathcal{G}_\mathbf{p},\mathcal{G}_\mathbf{q})$ which
is regular with respect to a Radon measure $\mathfrak{m}$ on
$\real_+$. Then there exist three unique families of adapted
processes $\{E_{ij}\},$ $\{F_i^*\},$ $\{G_i\}$ in
$\mathcal{L}(\mathcal{G}_\mathbf{p},\mathcal{G}_\mathbf{q})$ such
that
\begin{equation*}
{\rm d}\Xi =\sum_{i,j}E_{ij}{\rm d}\Lambda_{ij}+\sum_{i}F_i^*{\rm
d}A_i+\sum_{i}G_{i}{\rm d}A_i^*
\end{equation*}
on $\widetilde{\mathcal{E}}_{00}$ {\rm (}see eq.~$\eqref{eqn:tilde
of E00}${\rm )}. Furthermore{\rm ,}
\begin{equation*}
\sum_{i}\rho_i^{2q_3}G_i^\dagger(s)\mathbf{A}^{2\mathbf{q}}G_i(s)
\quad\text{and}\quad
\sum_{i}\rho_i^{-2p_3}F_i^\dagger(s)\mathbf{A}^{-2\mathbf{p}}F_i(s)
\end{equation*}
converge weakly to operators
$G(s)\in\mathcal{L}(\mathcal{G}_\mathbf{p},
\mathcal{G}_{-\mathbf{p}})$ and
$F(s)\in\mathcal{L}(\mathcal{G}_{-\mathbf{q}},
\mathcal{G}_{\mathbf{q}}),$ respectively{\rm ,} with
\begin{equation*}
\max\{\|G(s)\|_{\mathbf{p};-\mathbf{p}},\|F(s)\|_{-\mathbf{q};\mathbf{q}}\}
\le\mathfrak{m}_{\rm ac}'(s),\qquad s\in\real_+,
\end{equation*}
where $\mathfrak{m}_{\rm ac}$ denotes the absolutely continuous
part of $\mathfrak{m}$.
\end{theor}

\begin{proof}
This follows from the identity \eqref{eqn:IRs of annilation and
creation parts} and Lemma~\ref{lmm:integral representation against
Lambdaij} below.\hfill $\blacksquare$
\end{proof}

\begin{lem}\label{lmm:integrands against with Ai and Ai-dagger}
Let $\mathbf{p}=(p_1,p_2,p_3)$ and $\mathbf{q}=(q_1,q_2,q_3)$.
Let $\Xi$ be a martingale in
$\mathcal{L}(\mathcal{G}_\mathbf{p},\mathcal{G}_\mathbf{q})$
which is regular with respect to a Radon
measure $\mathfrak{m}$ on $\real_+$. Then
\begin{enumerate}
\renewcommand\labelenumi{\rm (\roman{enumi})}
\leftskip .4pc
\item $\mathfrak{m}$ can be replaced by
its absolutely continuous part{\rm ;}
\item there exist two countable
families of adapted processes $\{F_i^*(t)\}$ and $\{G_i(t)\}$ in
$\mathcal{L}(\mathcal{G}_\mathbf{p},\mathcal{G}_\mathbf{q})$ such
that for any $\varphi\in\mathcal{G}_{\mathbf{p};a]}$ and
$\psi\in\mathcal{G}_{-\mathbf{q};a]},$ $t>a\ge0,$
\begin{align*}
\hskip -1.25pc (\Xi(t)-\Xi(a))\varphi&= \int_a^t\,\sum_{i}G_i(s)\varphi\,{\rm d}B_i(s),\\[.3pc]
\hskip -1.25pc (\Xi^*(t)-\Xi^*(a))\psi&=
\int_a^t\,\sum_{i}F_i(s)\psi\,{\rm d}B_i(s),
\end{align*}
where $\{B_i(s)\}$ is the countable family of Brownian motions in
Proposition~$\ref{prop:extension of K-W Theorem};$
\item the series
\begin{equation*}
\hskip -1.25pc
\sum_{i}\rho_i^{2q_3}G_i^\dagger(s)\mathbf{A}^{2\mathbf{q}}G_i(s)
\quad\text{and}\quad
\sum_{i}\rho_i^{-2p_3}F_i^\dagger(s)\mathbf{A}^{-2\mathbf{p}}F_i(s)
\end{equation*}
converge weakly to operators
$G(s)\in\mathcal{L}(\mathcal{G}_\mathbf{p},\mathcal{G}_{-\mathbf{p}})$
and
$F(s)\in\mathcal{L}(\mathcal{G}_{-\mathbf{q}},\mathcal{G}_{\mathbf{q}}),$
respectively{\rm ,} with the property that
$\|G(s)\|_{\mathbf{p};-\mathbf{p}}$ and
$\|F(s)\|_{-\mathbf{q};\mathbf{q}}$ are locally integrable.
\end{enumerate}
\end{lem}

\begin{proof}$\left.\right.$

\noindent (i) Let $\varphi\in\mathcal{G}_{\mathbf{p};a]}$ be
fixed. Since $\{\mathbf{A}^{\mathbf{q}}\Xi(t)\varphi\}_{t\ge a}$
is a classical $\mathcal{I}$-valued square integrable martingale
in $[a,\infty)$ adapted to the countable family $\{B_i(s)\}$ of
independent Brownian motions, by Proposition \ref{prop:extension
of K-W Theorem} there exists a countable family of
$\mathcal{I}$-valued adapted square integrable process
$\{\xi_i(s,\varphi)\}_{s\ge a}$ such that
%
\begin{equation*}
\mathbf{A}^{\mathbf{q}}\Xi(t)\varphi-\mathbf{A}^{\mathbf{q}}\Xi(a)\varphi
=\int_a^t\sum_{i}\xi_i(s,\varphi){\rm d}B_i(s).
\end{equation*}
By the It\^o isometry and \eqref{eqn:regularity condition}, we
have for all $0\le a<b<t<\infty$,
\begin{align}\label{eqn:inequality for ac radon measure m}
\int_b^t\mathbf{E}\!\left[\sum_{i}\|\xi_i(s,\varphi)\|_\mathcal{I}^2\right]{\rm
d}s
=\pNORMM{0}{[\mathbf{A}^{\mathbf{q}}\Xi(t)\!-\!\mathbf{A}^{\mathbf{q}}\Xi(b)]\varphi}^2
\le\pNORM{\mathbf{p}}{\varphi}^2\mathfrak{m}([b,t]).
\end{align}
Similarly, we prove that for fixed
$\psi\in\mathcal{G}_{-\mathbf{q};a]}$ there exists a countable
family of $\mathcal{I}$-valued adapted square integrable process
$\{\eta_i(s,\psi)\}_{s\ge a}$ such that
\begin{equation*}
\mathbf{A}^{-\mathbf{p}}\Xi(t)^*\psi-\mathbf{A}^{-\mathbf{p}}\Xi(a)^*\psi
=\int_a^t\sum_{i}\eta_i(s,\psi){\rm d}B_i(s)
\end{equation*}
and for all $0\le a<b<t<\infty$,
\begin{align}\label{eqn:inequality for ac radon measure m2}
\int_b^t\mathbf{E}\left[\sum_{i}\|\eta_i(s,\psi)\|_\mathcal{I}^2\right]{\rm
d}s
&=\pNORMM{0}{[\mathbf{A}^{-\mathbf{p}}\Xi(t)^*-\mathbf{A}^{-\mathbf{p}}\Xi(b)^*]\psi}^2\nonumber\\[.3pc]
&\le\pNORM{-\mathbf{q}}{\psi}^2\mathfrak{m}([b,t]).
\end{align}
From \eqref{eqn:inequality for ac radon measure m} and
\eqref{eqn:inequality for ac radon measure m2} we see that
$\mathfrak{m}$ can be replaced by its absolutely continuous part
$\mathfrak{m}_{\rm ac}$.

\noindent (ii)--(iii). From (i) we assume that $\mathfrak{m}$ is
an absolutely continuous Radon measure. By similar arguments of
those used in the proof of Proposition~7.5 in \cite{Ji03} we see
that $\{\xi_i(s,\varphi)\}_{s\ge a}$ does not depend on the end
point $a$ and we put
\begin{equation*}
G_i(s)\varphi
=e^{-2q_2}\rho_i^{-q_3}\mathbf{A}^{-\mathbf{q}}\xi_i(s,\varphi)\quad\text{
a.e. }\quad s>a, \quad\varphi\in\mathcal{G}_{\mathbf{p};a]}.
\end{equation*}
This gives an adapted operator family $\{G_i(s)\}$ (for details
see the proof of the Proposition~7.5 in \cite{Ji03}). Hence by
\eqref{eqn:inequality for ac radon measure m} for any
$\varphi\in\mathcal{G}_{\mathbf{p};a]}$ we have
\begin{align*}
\int_b^t\sum_{i}\rho_i^{2q_3}\pNORMM{0}{\mathbf{A}^\mathbf{q}G_i(s)\varphi}^2{\rm
d}s
&=e^{-2q_2}\int_b^t\mathbf{E}\left[\sum_{i}\|\xi_i(s,\varphi)\|_\mathcal{I}^2\right]{\rm
d}s\\[.3pc] &\le
e^{-2q_2}\pNORM{\mathbf{p}}{\varphi}^2\mathfrak{m}([b,t])
\end{align*}
which implies that
\begin{align}\label{eqn:norm estimate of operator family Gi}
\sum_{i}\rho_i^{2q_3}\pNORMM{0}{\mathbf{A}^\mathbf{q}G_i(s)\varphi}^2
\le e^{-2q_2}\mathfrak{m}'(s)\pNORM{\mathbf{p}}{\varphi}^2
\quad\text{ for all }s.
\end{align}
This shows that each $G_i(s)$ is an adapted process in
$\mathcal{L}(\mathcal{G}_\mathbf{p},\mathcal{G}_\mathbf{q})$ and
that $\sum_i\rho_i^{2q_3}G_i^\dagger(s)\mathbf{A}^{2\mathbf{q}}$
$G_i(s)$ converges strongly to an operator
$G(s)\in\mathcal{L}(\mathcal{G}_\mathbf{p},\mathcal{G}_{-\mathbf{p}})$
such that $\|G(s)\|_{\mathbf{p};-\mathbf{p}}$ is locally
integrable. In fact, we prove that
\begin{equation*}
\left\|\sum_{i}
\rho_i^{2q_3}G_i^\dagger(s)\mathbf{A}^{2\mathbf{q}}G_i(s)
   \right\|_{\mathbf{p};-\mathbf{p}}
\le e^{-2q_2}\mathfrak{m}'(s)\quad\text{ a.e. }
\end{equation*}
The remainder of the proof is similar.\hfill $\blacksquare$
\end{proof}

Now, we put
\begin{align}
S(t)&=\int_0^t\sum_{i}\left(G_i(s){\rm d}A_i^*(s)+F_i^*(s){\rm
d}A_i(s)\right),\nonumber\\[.3pc]
S^*(t)&=\int_0^t\sum_{i}\left(F_i(s){\rm d}A_i^*(s)+G_i^*(s){\rm
d}A_i(s)\right),\nonumber\\[.3pc]
Z(t)&=\Xi(t)-S(t) \quad \textrm{ and }\quad
Z^*(t)=\Xi^*(t)-S^*(t).\label{eqn:IRs of annilation and creation
parts}
\end{align}

\begin{rem}\label{rem:integrals on tilde E0}
{\rm By \eqref{eqn:condition for QSI} and \eqref{eqn:norm estimate
of operator family Gi}, the integrals
\begin{equation*}
\int_0^t\sum_{i} G_i(s){\rm d}A_i^*, \quad
\int_0^t\sum_{i}F_i(s){\rm d}A_i^*
\end{equation*}
are well-defined on $\widetilde{\mathcal{E}}$ with $M=H$. But, in
general, the other two integrals $\int_0^t\sum_i F_i^*(s)$ ${\rm
d}A_i(s)$ and $\int_0^t\sum_i G_i^*(s){\rm d}A_i(s)$ are not
well-defined on $\widetilde{\mathcal{E}}$ with $M=H$ since we have
no estimates for
$\sum_{i}\rho_i^{-2p_3}\pNORMM{0}{\mathbf{A}^{-\mathbf{p}}G_i^*(s)\varphi}^2$
and
$\sum_{i}\rho_i^{2q_3}\pNORMM{0}{\mathbf{A}^{\mathbf{q}}F_i^*(s)\varphi}^2$.
If we consider the integrals on $\widetilde{\mathcal{E}}_0$, then
the infinite series reduce to finite sums and hence the stochastic
integrals are well-defined on $\widetilde{\mathcal{E}}_0$ by
\eqref{eqn:condition for QSI}. Then from \eqref{eqn:Bilinear form
of QSI} and the definitions it is immediate that the processes
$\{S,S^*\}$ and $\{Z,Z^*\}$ are adjoint pairs on
$\widetilde{\mathcal{E}}_0$. Also, we can easily see that for all
$u\in\mathcal{I}_\infty$ and $f\in M_0$,
$\{\mathbf{A}^\mathbf{q}Z(t)u\otimes\phi_{\mathbf{1}_{[0,t]}f}\}$
and
$\{\mathbf{A}^{-\mathbf{p}}Z^*(t)u\otimes\phi_{\mathbf{1}_{[0,t]}f}\}$
are classical $\mathcal{I}$-valued martingales adapted to the
countable family of Brownian motions $\{B_i\}$ in
Proposition~\ref{prop:extension of K-W Theorem}. Moreover, for all
$t>a$,
\begin{equation}\label{eqn:martingale property of Zt}
\hskip -4pc Z(t)u\otimes\phi_{\mathbf{1}_{[0,a]}f} =
Z(a)u\otimes\phi_{\mathbf{1}_{[0,a]}f},
Z^*(t)u\otimes\phi_{\mathbf{1}_{[0,a]}f} =
Z^*(a)u\otimes\phi_{\mathbf{1}_{[0,a]}f}.
\end{equation}}
\end{rem}

\begin{lem}\label{lmm:IR of Z(t)u-phi}
Let $u\in\mathcal{I}_\infty$ and $f\in M_0$. Then
\begin{enumerate}
\renewcommand\labelenumi{\rm (\roman{enumi})}
\leftskip .2pc
    \item there exists a $\mathcal{I}$-valued
                        square integrable classical process $\{\xi_i(\cdot,u,f)\}$
                        such that
\begin{align*}
\hskip -1.25pc
\mathbf{A}^\mathbf{q}Z(t)\mathbf{A}^{-\mathbf{p}}u\otimes\phi_{\mathbf{1}_{[0,t]}f}
=\mathbf{A}^\mathbf{q}\Xi(0)\mathbf{A}^{-\mathbf{p}}u\otimes\phi_0+\int_0^t
\sum_{i}\xi_i(s,u,f){\rm d}B_i(s);
\end{align*}

    \item there exists a $\mathcal{I}$-valued
                        square integrable classical process $\{\eta_i(\cdot,u,f)\}$
                        such that
\begin{align*}
\hskip -1.25pc
\mathbf{A}^{-\mathbf{p}}Z^*(t)\mathbf{A}^\mathbf{q}u\otimes\phi_{\mathbf{1}_{[0,t]}f}
=\mathbf{A}^{\!-\mathbf{p}}\Xi^*(0)\mathbf{A}^\mathbf{q}u\otimes\phi_0+\!\int_0^t
\sum_{i}\eta_i(s,u,f){\rm d}B_i(s).
\end{align*}
\end{enumerate}
\end{lem}

\begin{proof}
The proofs of (i) and (ii) are simple applications of
Proposition~\ref{prop:extension of K-W Theorem}.\hfill
$\blacksquare$
\end{proof}

Now, we prove that $\{Z\}_{t\ge0}$ can be represented by a
stochastic integral with respect to $\{\Lambda_{ij}\}$. For the
proof, we use similar arguments to those used in \cite{PS88} by
using a special martingale $U^{(i)}$ related to the Weyl
representation.

\begin{lem}\hskip -.2pc{\rm \cite{PS88}.}\ \ \label{lmm:martingale Ui}
For each $i=1,2,\ldots,$ let $U^{(i)}$ be the unique bounded
martingale satisfying
\begin{equation*}
\hbox{\rm d} U^{(i)}=({\rm d}A_i^*-{\rm d}A_i)U^{(i)}, \quad
U^{(i)}(0)=I.
\end{equation*}
Then
\begin{enumerate}
\renewcommand\labelenumi{\rm (\roman{enumi})}
\leftskip .2pc
\item $e^{-t/2}U^{(i)}(t)=I_0\otimes
W(\mathbf{1}_{[0,t]}e_i,I),$ where $I_0$ is the identity in
$\mathcal{B}(\mathcal{I})$ and $W$ is the Weyl representation
defined in {\rm \cite{HP84}}$;$

\item $U^{(i)}(t)$ leaves
$\widetilde{\mathcal{E}}_0$ invariant.
\end{enumerate}
\end{lem}

\begin{lem}\label{lmm:process Yi}
Let $\Xi$ be a regular martingale in
$\mathcal{L}(\mathcal{G}_\mathbf{p},\mathcal{G}_\mathbf{q})$ and
let $\{G_i\},$ $\{F_i\}$ be the associated families of adapted
processes defined in Lemma~$\ref{lmm:integrands against with Ai
and Ai-dagger}$. For each $i=1,2,\ldots,$ put
\begin{align*}
Y^{(i)}(t)
=\mathbf{A}^{\mathbf{q}}\left(\Xi(t)\mathbf{A}^{-\mathbf{p}}U^{(i)}(t)
    -e^{-p_2}\rho_i^{-p_3}\int_0^t F_i^*(s)\mathbf{A}^{-\mathbf{p}}U^{(i)}(s){\rm d}s\right).
\end{align*}
Then
\begin{enumerate}
\renewcommand\labelenumi{\rm (\roman{enumi})}
\leftskip .2pc
    \item for each $i,$ $\{Y^{(i)}\}_{t\ge0}$ is a bounded regular martingale
                      in $\mathcal{L}(\mathcal{G},\mathcal{G});$
    \item for each $i,$ there exists a unique family $\{M_j^{(i)}\}$
    of bounded adapted processes such that for all $t>a>0$ and
    $\varphi\in\mathcal{G}_{a]}$,
\begin{align}\label{eqn:IR of martingale generated by Yi}
\hskip -1.25pc [Y^{(i)}(t)- Y^{(i)}(a)]\varphi
=\int_a^t\sum_{j}M_j^{(i)}(s)\varphi {\rm d}B_j(s).
\end{align}
\end{enumerate}
\end{lem}

\begin{proof}$\left.\right.$

\noindent (i) It is clear that
$\mathbf{A}^\mathbf{q}\Xi(t)\mathbf{A}^{-\mathbf{p}}U^{(i)}(t)$ is
bounded. By similar arguments of those used to get \eqref{eqn:norm
estimate of operator family Gi} we prove that for any $t>0$,
\begin{equation}\label{eqn:norm estimate of Fis}
\int_0^t\|\mathbf{A}^{-\mathbf{p}}F_i(s)\mathbf{A}^{\mathbf{q}}U^{(i)}(s)\|_{0;0}{\rm
d}s \le e^{p_2}\int_0^t e^{{s}/{2}}\,\sqrt{\mathfrak{m}'(s)}{\rm
d}s<\infty
\end{equation}
which implies that $Y^{(i)}(t)$ is bounded. Since
$\widetilde{\mathcal{E}}_0$ is invariant by $U^{(i)}(t)$, the
relation:
\begin{equation}\label{eqn:another represen of Yi}
Y^{(i)}(t)
=\mathbf{A}^\mathbf{q}Z(t)\mathbf{A}^{-\mathbf{p}}U^{(i)}(t)+W^{(i)}(t)
\end{equation}
holds on $\widetilde{\mathcal{E}}_0$, where
\begin{align*}
&W^{(i)}(t)\\[.4pc]
&\quad\,=\mathbf{A}^{\mathbf{q}}\left(S(t)\mathbf{A}^{-\mathbf{p}}U^{(i)}(t)
-e^{-p_2}\rho_i^{-p_3}\int_0^tF_i^*(s)\mathbf{A}^{-\mathbf{p}}U^{(i)}(s){\rm
d}s\right), \quad t\ge0.
\end{align*}
Note that $\widetilde{\mathcal{E}}_0$ is invariant by
$\mathbf{A}^\mathbf{p}$. Therefore, by \eqref{eqn:martingale
property of Zt} and the martingale property of $U^{(i)}$, we prove
that
$\{\mathbf{A}^{\mathbf{q}}Z(t)\mathbf{A}^{-\mathbf{p}}U^{(i)}(t)\}$
is a martingale on $\widetilde{\mathcal{E}}_0$. Now let
$u,v\in\mathcal{I}_\infty$ and $f,g \in M_0$. Then by applying
Ito's product formula \eqref{eqn:quantum Ito formula} (or see
Theorem~6.2 in \cite{Ji03}) to
$\mathbf{A}^{\mathbf{q}}S(t)\mathbf{A}^{-\mathbf{p}}U^{(i)}(t)$,
we can easily see that $\{W^{(i)}(t)\}_{t\ge0}$ is a martingale.

In fact, for any $t\ge0$ we have
\begin{align}
\mathbf{A}^{\mathbf{q}}S(t)\mathbf{A}^{-\mathbf{p}}U^{(i)}(t)
 &=e^{-p_2}\rho_i^{-p_3}\int_0^t\mathbf{A}^{\mathbf{q}}S(s)\mathbf{A}^{-\mathbf{p}}U^{(i)}(s){\rm d}A_i^*(s)\nonumber\\[.3pc]
 &\quad\,-\int_0^t\mathbf{A}^{\mathbf{q}}S(s)\mathbf{A}^{-\mathbf{p}}U^{(i)}(s){\rm d}A_i(s)
      \nonumber\\[.3pc]
 &\quad\,+e^{q_2}\rho_i^{q_3}\int_0^t\mathbf{A}^{\mathbf{q}}G_i(s)\mathbf{A}^{-\mathbf{p}}U^{(i)}(s){\rm d}A_i^*(s)
       \nonumber\\[.3pc]
 &\quad\,+\int_0^t\mathbf{A}^{\mathbf{q}}F_i^*(s)\mathbf{A}^{-\mathbf{p}}U^{(i)}(s){\rm d}A_i(s)
       \nonumber\\[.3pc]
 &\quad\,+e^{-p_2}\rho_i^{-p_3}\int_0^t\mathbf{A}^{\mathbf{q}}F_i^*(s)\mathbf{A}^{-\mathbf{p}}U^{(i)}(s){\rm d}s
  \label{eqn:S(t)Ui(t)}
\end{align}
on $\widetilde{\mathcal{E}}_0$. The proof of regularity is similar
to that in \cite{PS86} and \cite{PS88}.  By similar arguments of
those used to get \eqref{eqn:norm estimate of Fis} we first show
that for $t>a>0$ and $\varphi\in\mathcal{G}_{a]}$,
\begin{align*}
&\pNORM{0}{e^{-p_2}\rho_i^{-p_3}\int_a^tF_i^*(s)\mathbf{A}^{-\mathbf{p}}U^{(i)}(s)\varphi
{\rm d}s}^2\\[.3pc]
&\quad\, \le \rho_i^{-2p_3}\pNORM{0}{\varphi}^2 \left(\int_a^t
e^{s/2}\sqrt{\mathfrak{m}'(s)}{\rm d}s\right)^2\\[.3pc]
&\quad\, \le
\rho_i^{-2p_3}\pNORM{0}{\varphi}^2(e^t-e^a)\mathfrak{m}([a,t]).
\end{align*}
On the other hand, for $t>a>0$ and $\varphi\in\mathcal{G}_{a]}$ we
have
\begin{align*}
&\pNORMM{q}{(\Xi(t)\mathbf{A}^{-\mathbf{p}}U^{(i)}(t)
-\Xi(a)\mathbf{A}^{-\mathbf{p}}U^{(i)}(a))\varphi}^2\\[.3pc]
&\quad\, \le 2\pNORMM{q}{\Xi(t)\mathbf{A}^{-\mathbf{p}}(U^{(i)}(t)
  -U^{(i)}(a))\varphi}^2\\[.3pc]
&\qquad\, +2\pNORMM{q}{(\Xi(t)-\Xi(a))\mathbf{A}^{-\mathbf{p}}U^{(i)}(a)\varphi}^2\\[.3pc]
&\quad\, \le 2\|\Xi(t)\|_{\mathbf{p};\mathbf{q}}^2
(e^t-e^a)\pNORM{0}{\varphi}^2
+2e^a\pNORM{0}{\varphi}^2\mathfrak{m}([a,t]).
\end{align*}
Therefore, since $\|\Xi(t)\|_{\mathbf{p};\mathbf{q}}$ is
non-decreasing by Remark~\ref{remark:remark1 for IRT}, for $t>a>0$
and $\varphi\in\mathcal{G}_{a]}$ we have
\begin{align}
&\pNORMM{0}{(Y^{(i)}(t) - Y^{(i)}(a))\varphi}^2\nonumber\\[.3pc]
&\quad\,\le2C\pNORM{0}{\varphi}^2((e^t-e^a)\mathfrak{m}([a,t])
  +(e^t-e^a)\|\Xi(t)\|_{\mathbf{p};\mathbf{q}}^2
  +e^a\mathfrak{m}([a,t]))\nonumber\\[.3pc]
&\quad\,\le2C\pNORM{0}{\varphi}^2[e^t(\|\Xi(t)\|_{\mathbf{p};\mathbf{q}}^2
            +\mathfrak{m}([0,t]))
-e^a(\|\Xi(a)\|_{\mathbf{p};\mathbf{q}}^2+\mathfrak{m}([0,a]))]\nonumber\\[.3pc]
&\quad\,\equiv \pNORM{0}{\varphi}^2\mathfrak{n}([a,t]),
\label{eqn:regularity of Yi}
\end{align}
where $C=\max\{\rho_i^{-2p_3},2\}$ and $\mathfrak{n}$ is the Radon
measure defined by \eqref{eqn:regularity of Yi}. Hence we prove
(i). The proof of (ii) is similar to the proof of
Lemma~\ref{lmm:integrands against with Ai and Ai-dagger} by
applying Proposition~ \ref{prop:extension of K-W Theorem} to the
bounded regular martingale $\{Y^{(i)}(t)\}_{t\ge0}$.\hfill
$\blacksquare$
\end{proof}

Let
\begin{equation*}
M_{00}=\{f\in M_0\subset H_\infty\,;\, \|f(t)\|\text{ is a locally
bounded function of }t\}
\end{equation*}
and
\begin{equation}\label{eqn:tilde of E00}
\widetilde{\mathcal{E}}_{00}=\mathcal{I}_\infty\otimes_{\rm
al}\mathcal{E}(M_{00}).
\end{equation}

\begin{lem}\label{lmm:integral representation against Lambdaij}
Let $\Xi$ be a regular martingale in
$\mathcal{L}(\mathcal{G}_\mathbf{p},\mathcal{G}_\mathbf{q})$ and
let $\{F_i^*\},$ $\{G_i\},$ $\{S,S^*\},$ $\{Z,Z^*\},$ $U^{(i)}$
and $\{M_j^{(i)}\}$ be as defined in Lemmas~$\ref{lmm:integrands
against with Ai and Ai-dagger}$--$\ref{lmm:process Yi}$. Put
\begin{align}\label{eqn:definition of E-ij}
L_{ij}(t)&=
\mathbf{A}^{-\mathbf{q}}M_j^{(i)}(t)U^{(i)}(t)^{-1}\mathbf{A}^{\mathbf{p}}\nonumber\\[.3pc]
&\quad\,-e^{q_2}\rho_i^{q_3}G_i(t)-e^{-p_2}\rho_i^{-p_3}S(t)\delta_{ij}-Z(t)\delta_{ij}.
\end{align}
Then the processes $\{\eta_i(t,u,f)\}$ defined in {\rm (ii)} of
Lemma~$\ref{lmm:IR of Z(t)u-phi}$ satisfies the relation for any
$u \in\mathcal{I}_\infty$ and $f\in M_0${\rm :}
\begin{equation}\label{eqn:eta-i(t,u,f)}
\eta_i(t,u,f) =\sum_{j}
  f_j(t)\mathbf{A}^{-\mathbf{p}}[L_{ij}^*(t) + Z^*(t)\delta_{ij}]\mathbf{A}^{\mathbf{q}}
  u\otimes\phi_{\mathbf{1}_{[0,t]}f}\quad\text{ a.e. }\,\,t.
\end{equation}
Moreover{\rm ,} we have
\begin{equation}\label{eqn:Z(t)}
 Z(t)=\Xi(0)+\int_0^t\sum_{i,j}E_{ij}(s){\rm d}\Lambda_{ij}(s)
\end{equation}
defined on $\widetilde{\mathcal{E}}_{00},$ where
$E_{ij}=e^{p_2-q_2}\rho_i^{-q_3}\rho_j^{p_3}L_{ij}$ for each
$i,j=1,2,\ldots$.
\end{lem}

\begin{proof}
Let $u,v\in\mathcal{I}_\infty$, $f,g \in M_0$ and $t>a$. Then by
(i) in Lemma \ref{lmm:martingale Ui} we have
\begin{equation*}
U^{(i)}(t)v\otimes\phi_{\mathbf{1}_{[0,a]}g} ={\rm
e}^{-\int_0^ag_i(s){\rm d}s}
v\otimes\phi_{\mathbf{1}_{[0,a]}g+\mathbf{1}_{[0,t]}e_i}.
\end{equation*}
Thus by (ii) in Lemma~\ref{lmm:IR of Z(t)u-phi} and the It\^o
isometry, we have
\begin{align*}
\hskip -4pc&\frac{\rm d}{{\rm
d}t}\Bilinn{\mathbf{A}^{-\mathbf{p}}Z^*(t)\mathbf{A}^\mathbf{q}u\otimes\phi_{\mathbf{1}_{[0,t]}f}}
    {U^{(i)}(t)v\otimes\phi_{\mathbf{1}_{[0,a]}g}}\nonumber\\[.3pc]
\hskip -4pc &\quad\, =\hbox{e}^{-\int_0^ag_i(s){\rm d}s} \frac{\rm
d}{{\rm d}t}\int_0^t\sum_{j}
 (\mathbf{1}_{[0,a]}(s)g_i(s)+ \delta_{ij})
 \Bilinn{\eta_j(s,u,f)}{v\otimes\phi_{\mathbf{1}_{[0,a]}g+\mathbf{1}_{[0,t]}e_i}}{\rm d}s\nonumber\\[.3pc]
\hskip -4pc &\quad\,
=\Bilinn{\eta_i(t,u,f)}{U^{(i)}(t)v\otimes\phi_{\mathbf{1}_{[0,a]}g}}
 \quad\text{a.e.}\,\,t>a.\\[-4.1pc]
\hskip -4pc
\end{align*}
\begin{align}
\label{eqn:ZtUit}
\end{align}\vspace{.4pc}

\noindent On the other hand, from \eqref{eqn:another represen of
Yi}, \eqref{eqn:S(t)Ui(t)}, \eqref{eqn:IR of martingale generated
by Yi} and It\^o isometry, for $t > a$ we obtain that
\begin{align*}
\hskip -4pc &\Bilinn{u\otimes\phi_{\mathbf{1}_{[0,t]}f}}
   {[\mathbf{A}^\mathbf{q}Z(t)\mathbf{A}^{-\mathbf{p}}U^{(i)}(t)
                 -\mathbf{A}^\mathbf{q}Z(a)\mathbf{A}^{-\mathbf{p}}U^{(i)}(a)]
          v\otimes\phi_{\mathbf{1}_{[0,a]}g}}\nonumber\\[.3pc]
\hskip -4pc &\quad\,=\Bilinn{u\otimes\phi_{\mathbf{1}_{[0,t]}f}}
   {[Y^{(i)}(t)-Y^{(i)}(a)]v\otimes\phi_{\mathbf{1}_{[0,a]}g}}\nonumber\\[.3pc]
\hskip -4pc &\qquad\, -\Bilinn{u\otimes\phi_{\mathbf{1}_{[0,t]}f}}
   {[W^{(i)}(t)-W^{(i)}(a)]v\otimes\phi_{\mathbf{1}_{[0,a]}g}}\nonumber
\end{align*}\pagebreak
\begin{align*}
\hskip -4pc &\quad\,
  =\int_a^t \sum_{j} f_j(s)
 \Bilinn{u\otimes\phi_{\mathbf{1}_{[0,t]}f}}{M_j^{(i)}(s)v\otimes\phi_{\mathbf{1}_{[0,a]}g}}{\rm d}s
    -\int_a^t\sum_{j}{f}_j(s)\nonumber\\[.3pc]
\hskip -4pc &\qquad\,\times
    \Bilinn{u\otimes\phi_{\mathbf{1}_{[0,t]}f}}
    {\mathbf{A}^\mathbf{q}[e^{q_2}\rho_i^{q_3}G_i(s)+e^{-p_2}\rho_i^{-p_3}S(s)\delta_{ij}]
       \mathbf{A}^{-\mathbf{p}}U^{(i)}(s)v\otimes\phi_{\mathbf{1}_{[0,a]}g}}{\rm d}s\nonumber\\[.3pc]
\hskip -4pc &\quad\,
  =\int_a^t\sum_{j}{f}_j(s)
   \Bilinn{u\otimes\phi_{\mathbf{1}_{[0,t]}f}}
   {\mathbf{A}^\mathbf{q}[L_{ij}(s)+Z(s)\delta_{ij}]\mathbf{A}^{-\mathbf{p}}
       U^{(i)}(s)v\otimes\phi_{\mathbf{1}_{[0,a]}g}}{\rm d}s.\\[-3pc]
\hskip -4pc
\end{align*}
\begin{align}
\label{eqn:ZtUit-2}
\end{align} \noindent Therefore, by comparing \eqref{eqn:ZtUit} and
\eqref{eqn:ZtUit-2} using the totality of the set
$\{v\otimes\phi_{\mathbf{1}_{[0,a]}g}\,|\,v\in\mathcal{I}_\infty,\,\,g\in{M}_0,\,\,0<a<t\}$
in $\mathcal{G}_{t]}$, we have \eqref{eqn:eta-i(t,u,f)}. Hence by
(ii) in Lemma \ref{lmm:IR of Z(t)u-phi} and
 \eqref{eqn:eta-i(t,u,f)} we obtain that
\begin{align*}
\hskip -4pc &\mathbf{A}^{-\mathbf{p}}Z^*(t)\mathbf{A}^\mathbf{q}u\otimes\phi_{\mathbf{1}_{[0,t]}f}\nonumber\\[.3pc]
\hskip -4pc &\quad
=\mathbf{A}^{-\mathbf{p}}\Xi^*(0)\mathbf{A}^\mathbf{q}u\otimes\phi_0
+\int_0^t\sum_{i,j}f_j(s)
\mathbf{A}^{-\mathbf{p}}[L_{ij}^*(s)+Z^*(s)\delta_{ij}]\mathbf{A}^{\mathbf{q}}
u\otimes\phi_{\mathbf{1}_{[0,s]}f}{\rm d}B_i(s).\\[-3pc]
\hskip -4pc
\end{align*}
\begin{align}
\label{eqn:Z*(t)}
\end{align}
It is obvious that $\{L_{ij}(t)\}_{t\ge0}$ defined by
\eqref{eqn:definition of E-ij} are adapted processes in
$\mathcal{L}(\mathcal{G}_\mathbf{p},\mathcal{G}_\mathbf{q})$ and a
simple estimate shows that the integral
$\int_0^t\sum_{i,j}L_{ij}(s){\rm d}\Lambda_{ij}(s)$ is
well-defined on $\widetilde{\mathcal{E}}_{00}$ since the
integrability condition \eqref{eqn:condition for QSI} is satisfied
for any $f\in{M}_{00}$. Now, by \eqref{eqn:Z*(t)} and the It\^o
isometry, for $f,g\in{M}_{00}$ we have
\begin{align*}
\hskip -4pc &\Bilinn{u\otimes\phi_f}{\mathbf{A}^\mathbf{q}Z(t)\mathbf{A}^{-\mathbf{p}}v\otimes\phi_g}\\[.4pc]
\hskip -4pc &\quad\,
  =\Bilinn{\mathbf{A}^{-\mathbf{p}}Z^*(t)\mathbf{A}^\mathbf{q}u\otimes\phi_{\mathbf{1}_{[0,t]}f}}
         {v\otimes\phi_{\mathbf{1}_{[0,t]}g}}{\rm e}^{\int_t^\infty\bilin{f(s)}{g(s)}{\rm d}s}\\[.4pc]
\hskip -4pc &\quad\,
  ={\rm e}^{\int_t^\infty\bilin{f(s)}{g(s)}{\rm d}s}
   \Bigg\{\Bilinn{u\otimes\phi_0}{\mathbf{A}^\mathbf{q}\Xi(0)\mathbf{A}^{-\mathbf{p}}v\otimes\phi_0}\\[.3pc]
\hskip -4pc &\qquad\,
    +\int_0^t\sum_{i,j}{f}_j(s)g_i(s)
   \Bilinn{u\otimes\phi_{\mathbf{1}_{[0,s]}f}}
      {\mathbf{A}^{\mathbf{q}}[L_{ij}(s)+Z(s)\delta_{ij}]\mathbf{A}^{-\mathbf{p}}
        v\otimes\phi_{\mathbf{1}_{[0,s]}g}}{\rm d}s\Bigg\}.
\end{align*}
By differentiation we obtain that
\begin{align*}
&\frac{{\rm d}}{{\rm
d}t}\Bilinn{u\otimes\phi_f}{\mathbf{A}^\mathbf{q}Z(t)\mathbf{A}^{-\mathbf{p}}v\otimes\phi_g}\\[.3pc]
&\quad\, =\sum_{i,j}{f}_j(t)g_i(t)
 \Bilinn{u\otimes\phi_f}{\mathbf{A}^\mathbf{q}L_{ij}(t)\mathbf{A}^{-\mathbf{p}}v\otimes\phi_g}
\end{align*}
which by \eqref{eqn:Bilinear form of QSI} and (ii) in
Lemma~\ref{lmm:IR of Z(t)u-phi}, proves that
\begin{equation}\label{eqn:Z(t)-2}
Z(t)=\Xi(0)
+\mathbf{A}^{-\mathbf{q}}\left[\int_0^t\sum_{i,j}\mathbf{A}^{\mathbf{q}}L_{ij}(s)
\mathbf{A}^{-\mathbf{p}}{\rm
d}\Lambda_{ij}(s)\right]\mathbf{A}^{\mathbf{p}}
\end{equation}\pagebreak

\noindent on $\widetilde{\mathcal{E}}_{00}$. On the other hand, by
direct computation we prove that
\begin{align*}
&\mathbf{A}^{-\mathbf{q}}\left[\int_0^t\sum_{i,j}\mathbf{A}^{\mathbf{q}}L_{ij}(s)
\mathbf{A}^{-\mathbf{p}}{\rm
d}\Lambda_{ij}(s)\right]\mathbf{A}^{\mathbf{p}}\\[.4pc]
&\quad\, =\int_0^t\sum_{i,j} {\rm
e}^{p_2-q_2}\rho_i^{-q_3}\rho_j^{p_3}L_{ij}(s){\rm
d}\Lambda_{ij}(s)
\end{align*}
on $\widetilde{\mathcal{E}}_{00}$. Thus by \eqref{eqn:Z(t)-2} we
prove \eqref{eqn:Z(t)}.\hfill $\blacksquare$
\end{proof}

\section*{Acknowledgement}

This work was supported by the Korea Research Foundation Grant
funded by the Korean Government (MOEHRD) (No.
R05-2004-000-11346-0).

\end{document}